\documentclass[12pt]{article}

\usepackage{amsmath,amssymb,amsbsy,amsfonts,amsthm,latexsym,
                  amsopn,amstext,amsxtra,euscript,amscd}

\begin{document}

\newtheorem{theorem}{Theorem}
\newtheorem{lemma}[theorem]{Lemma}
\newtheorem{claim}[theorem]{Claim}
\newtheorem{cor}[theorem]{Corollary}
\newtheorem{prop}[theorem]{Proposition}
\newtheorem{definition}{Definition}
\newtheorem{question}[theorem]{Open Question}

\def\cA{{\mathcal A}}
\def\cB{{\mathcal B}}
\def\cC{{\mathcal C}}
\def\cD{{\mathcal D}}
\def\cE{{\mathcal E}}
\def\cF{{\mathcal F}}
\def\cG{{\mathcal G}}
\def\cH{{\mathcal H}}
\def\cI{{\mathcal I}}
\def\cJ{{\mathcal J}}
\def\cK{{\mathcal K}}
\def\cL{{\mathcal L}}
\def\cM{{\mathcal M}}
\def\cN{{\mathcal N}}
\def\cO{{\mathcal O}}
\def\cP{{\mathcal P}}
\def\cQ{{\mathcal Q}}
\def\cR{{\mathcal R}}
\def\cS{{\mathcal S}}
\def\cT{{\mathcal T}}
\def\cU{{\mathcal U}}
\def\cV{{\mathcal V}}
\def\cW{{\mathcal W}}
\def\cX{{\mathcal X}}
\def\cY{{\mathcal Y}}
\def\cZ{{\mathcal Z}}

\def\A{{\mathbb A}}
\def\B{{\mathbb B}}
\def\C{{\mathbb C}}
\def\D{{\mathbb D}}
\def\E{{\mathbb E}}
\def\F{{\mathbb F}}
\def\G{{\mathbb G}}
\def\H{{\mathbb H}}
\def\I{{\mathbb I}}
\def\J{{\mathbb J}}
\def\K{{\mathbb K}}
\def\L{{\mathbb L}}
\def\M{{\mathbb M}}
\def\N{{\mathbb N}}
\def\O{{\mathbb O}}
\def\P{{\mathbb P}}
\def\Q{{\mathbb Q}}
\def\R{{\mathbb R}}
\def\S{{\mathbb S}}
\def\T{{\mathbb T}}
\def\U{{\mathbb U}}
\def\V{{\mathbb V}}
\def\W{{\mathbb W}}
\def\X{{\mathbb X}}
\def\Y{{\mathbb Y}}
\def\Z{{\mathbb Z}}

\def\E{{\mathbf E}}
\def\Fp{\F_p}
\def\ep{{\mathbf{e}}_p}
\def\Bn{{\mathfrak B}_n}

\def\scr{\scriptstyle}
\def\\{\cr}
\def\({\left(}
\def\){\right)}
\def\[{\left[}
\def\]{\right]}
\def\<{\langle}
\def\>{\rangle}
\def\fl#1{\left\lfloor#1\right\rfloor}
\def\rf#1{\left\lceil#1\right\rceil}
\def\le{\leqslant}
\def\ge{\geqslant}
\def\eps{\varepsilon}
\def\mand{\qquad\mbox{and}\qquad}

\def\vec#1{\mathbf{#1}}
\def\inv#1{\overline{#1}}
\def\vol#1{\mathrm{vol}\,{#1}}

\newcommand{\comm}[1]{\marginpar{%
\vskip-\baselineskip 
\raggedright\footnotesize
\itshape\hrule\smallskip#1\par\smallskip\hrule}}

\def\xxx{\vskip5pt\hrule\vskip5pt}


\title{\bf On the Distribution of  Kloosterman 
Sums}

\author{ 
{\sc Igor E. Shparlinski} \\
{Department of Computing, Macquarie University} \\
{Sydney, NSW 2109, Australia} \\
{igor@ics.mq.edu.au}}

\date{\today}
\pagenumbering{arabic}

\maketitle

\begin{abstract}
For   a   prime $p$, we consider Kloosterman 
sums
$$
K_{p}(a) = \sum_{x\in \F_p^*} \exp( 2 \pi i (x + ax^{-1})/p), \qquad a \in \F_p^*,  
$$
over a finite field of $p$ elements.
It is well known that  due to results of Deligne, Katz and Sarnak, 
the distribution of the sums $K_{p}(a)$ when $a$ runs through $\F_p^*$
is in accordance with the Sato--Tate conjecture. 
Here we show that the same holds where $a$ runs through the sums
$a = u+v$ for $u \in \cU$, $v \in \cV$ for any two sufficiently large
sets $\cU, \cV \subseteq \F_p^*$. 

We also improve  a recent bound on the 
nonlinearity of a Boolean function associated with the sequence of
signs of Kloosterman sums. 
\end{abstract}

\section{Introduction}
\label{sec:intro}

For a prime $p$ we use $\F_p$ to denote the finite field of $p$ 
elements.  

For $a \in \F_p^*$ we consider the Kloosterman 
sum
$$
K_{p}(a) = \sum_{x\in \F_p^*} \ep(x + ax^{-1}), 
$$
where
$$
\ep(z) = \exp( 2 \pi i z/p)
$$
(we identify  $\F_p$ with the set $\{0, 1, \ldots, p-1\}$).
Since for the complex conjugated sum we have
$$
\inv{ K_{p}(a)}  = \sum_{x\in \F_p^*} \ep(-x - ax^{-1}) =  K_{p}(a) 
$$
the values of $K_{p}(a)$ are real. 

Accordingly to the Weil bound, see~\cite{LN}, 
$$
\left| K_{p}(a)\right| \le 2 \sqrt{p}, \qquad a \in \F_p^*. 
$$
Therefore,  we can  define the angles $\psi_{p}(a)$ 
by the relations
$$
 K_{p}(a) = 2\sqrt{p} \cos \psi_{p}(a) \mand 0 \le \psi_{p}(a)  \le \pi.
$$

The famous \emph{Sato--Tate} conjecture asserts that 
for any fixed non-zero integer $a$, when $p$ varies,  the angles
$\psi_{p}(a)$ are distributed accordingly to  the
\emph{Sato--Tate  density}
$$
 \mu_{ST}(\alpha,\beta) = \frac{2}{\pi}\int_\alpha^\beta  \sin^2 \gamma\, d \gamma, 
$$
see~\cite{Adolph,ChaLi,FoMic1,FoMic2,FMRS,Katz,KatzSar,Lau,Mich1,Mich2,Nied}
for various 
modifications and generalisations of this conjecture and further references.

It is also known that when a sufficiently large prime $p$ is fixed  
and $a$ runs through $\F_p^*$, then, as has been shown by Katz~\cite[Chapter~13]{Katz}, 
the work of Deligne on the \emph{Weil conjecture}  implies that the  
distribution of  the sums $K_{p}(a)$ is in accordance with the 
Sato--Tate density, 
see also~\cite{KatzSar}. Furthermore, a quantitative form of 
this result is given by Niederreiter~\cite{Nied}. 
Namely, if  $\cA_p(\alpha,\beta)$ is the set  
of $a \in \F_p^*$  with $ \alpha \le \psi_{p}(a) \le \beta$ 
then by the main result of Niederreiter~\cite{Nied}, we have: 
\begin{equation}
\label{eq:Nied Bound}
\max_{0 \le \alpha < \beta \le \pi} 
\left|\# \cA_p(\alpha,\beta) - \mu_{ST}(\alpha,\beta) p  \right| \ll p^{3/4}.
\end{equation}

Combining  results of  Fouvry, Michel, Rivat and S{\'a}rk{\"o}zy~\cite[Lemma~2.3]{FMRS} 
(with $r=1$)  and of Niederreiter~\cite[Lemma~3]{Nied}  one can show that elements of
$\cA_p(\alpha,\beta) $ are uniformly distributed, 
in the following sense. For any $\lambda \in \F_p^*$
and  integer $M$ with  $1 \le  M \le p-1$
we put 
$$
\cA_p(\lambda, M;\alpha,\beta) = \{a \in \cA_p(\alpha,\beta) \ : \ \lambda a \in [1,M]\}. 
$$
Then for $1 \le M \le p-1$, the following bound holds:
\begin{equation}
\label{eq:FMRS-Nied Bound}
\max_{\lambda \in \F_p^*}\, \max_{0 \le \alpha < \beta \le \pi} 
\left|\# \cA_p(\lambda, M;\alpha,\beta)  - \mu(\alpha,\beta) M  \right|
 \ll M^{1/2} p^{1/4} (\log p)^{1/2} .
\end{equation}

Fouvry, Michel, Rivat and S{\'a}rk{\"o}zy~\cite{FMRS} also remark that
combining a result of Fouvry and Michel~\cite{FoMic1}  with the technique of 
Vaaler~\cite{Vaal} one can show that
$$
 \max_{0 \le \alpha < \beta \le \pi} 
\left|\# \cQ_p(\alpha,\beta)  - \mu_{ST}(\alpha,\beta) p  \right| \ll p^{3/4} 
$$
where
$$
\cQ_p(\alpha,\beta)   = \{a \in \F_p : \ a^2 \in \cA_p(\alpha,\beta)\}. 
$$
The same bound can also be immediately obtained if one applies  the 
result of Niederreiter~\cite[Lemma~3]{Nied} to the bound  
of  Michel~\cite[Corollary~2.4]{Mich2} 
(see also~\cite[Lemma~2.1]{FoMic1}).

Here we show that the same type of distribution is preserved when $a$ runs through the sums
$a = u + v$  where $u \in \cU$, $v \in \cV$ for any two sufficiently large
sets $\cU, \cV \subseteq \F_p^*$. 
Namely, for any two sets $\cU, \cV \subseteq \F_p^*$, we put 
$$
\cW_p(\cU, \cV;\alpha,\beta) = \{(u,v) \in \cU\times \cV\ : \ u+v
\in \cA_p(\alpha,\beta)\}. 
$$
In particular, we obtain an asymptotic formula for 
$\#\cW_p(\cU, \cV;\alpha,\beta) $ which is nontrivial whenever 
\begin{equation}
\label{eq:Threshold}
\# \cU \# \cV \ge p^{3/2 + \varepsilon}
\end{equation} 
for any fixed $\varepsilon > 0$ and sufficiently large $p$. 

Then, we also improve the upper bound of~\cite{Shp} on the
\emph{nonlinearity} of the Boolean function
associated with the sequence  of signs of Kloosterman sums, that is
for the function
\begin{equation}
\label{eq:Bool}
f(a) = \left\{  \begin{array}{ll}
0,& \quad \text{if}\  K(a) > 0\ \text{or}\ a=0 , \\
1, & \quad \mbox{if}\ K(a) < 0,
\end{array} \right. \qquad a = 0, 1, \ldots\,2^{n-1}, 
\end{equation}
where $n$ is defined  by the inequalities
$$
2^n \le  p < 2^{n+1}.
$$

We denote by $\Bn$ the $n$-dimensional 
Boolean cube $\Bn = \{0,1\}^n$  and in a natural 
way identify 
its elements with the integers in the range $0 \le a \le 2^n -1$
(and thus with a subset of $\F_p$). 

We define the {\it Fourier
coefficients\/} of $f(a)$ as 
$$
{\widehat f}(r) = 2^{-n} \sum_{a \in\Bn}
(-1)^{  f(a) + \langle h, r  \rangle}  , \qquad r \in \Bn, 
$$
where   $\langle a,  r \rangle$ denote the inner product of
$a,r \in\Bn$. Furthermore, we recall that 
$$
N(f) =  2^{n-1} -  2^{n-1} \max_{r \in \Bn}\left| {\widehat f}(r)  \right|
$$
is called the {\it nonlinearity\/} of $f$ and is an important
cryptographic  characteristic, for example, see~\cite{CarDing}. 
In particular, it is   the smallest 
possible Hamming distance between the vector of values of $f$ and the vector 
of values of a linear function in $n$ variables over the $\F_2$.

Several results about some  measures of pseudorandomess 
of the sequence  of signs of Kloosterman sums have recently 
been obtained by 
Fouvry, Michel, Rivat and S{\'a}rk{\"o}zy~\cite{FMRS}.
Motivated by (and actually  using) the results of~\cite{FMRS}, 
the bound 
$$
N(f) = 2^{n-1} \(1  + O\( 2^{-n/24} n^{1/12} \)\).
$$
has been obtained in~\cite{Shp}.  
Here we again use some results of~\cite{FMRS}, 
but in a slightly different way and  improve this bound.

%
%

\section{Distribution of Elements of $\cA_p(\alpha,\beta)$}
\label{sec:prep}

For a sequence of
$N$ real numbers $\gamma_1,\ldots,\gamma_N\in[0,1)$ the  \emph{discrepancy} is
defined by
$$
D =\max_{0\le \gamma  \le 1}\left|\frac{T(\gamma,N)}{N}-\gamma \right|,
$$
where  $T(\gamma,N)$ is the number of  
$n\le N$ such that $\gamma_n\le \gamma$.

We also recall our agreements that the elements of $\F_p$ have 
canonical representation as integers of the interval $[0,p-1]$.
Thus for any field element $c \in\F_p$ we interpret $c/p$ as
a rational number in the interval $[0,1]$. 
Hence, for $\lambda \in \F_p^*$ we can define
the discrepancy $D_p(\lambda;\alpha,\beta)$ of the sequence
$$
\frac{\lambda a}{p}, \quad  a \in \cA_p(\alpha,\beta). 
$$ 
Then  the bound~\eqref{eq:FMRS-Nied Bound} implies that
$$
\max_{1 \le M \le p-1}\, 
\max_{\lambda \in \F_p^*}\, \max_{0 \le \alpha < \beta \le \pi} 
\left|\# \cA_p(\lambda, M;\alpha,\beta)  - \mu(\alpha,\beta) M  \right|
 \ll p^{3/4} (\log p)^{1/2}  
$$
which can be reformulated 
in the following form:

\begin{lemma}
\label{lem:Discr} We have,
$$
 \max_{\lambda \in \F_p^*}\,\max_{0 \le \alpha < \beta \le \pi} D_p(\lambda;\alpha,\beta) \ll 
p^{-1/4} (\log p)^{1/2} .
$$
\end{lemma}

Our main tool is a bound of exponential sums
 with elements
of $\cA_p(\alpha,\beta)$. For $\lambda \in \F_p^*$ we define
$$
S_p(\lambda;\alpha,\beta) = \sum_{a \in \cA_p(\alpha,\beta)} \ep(\lambda a).
$$

\begin{lemma}
\label{lem:Exp Sum} We have,
$$
 \max_{\lambda \in \F_p^*}\,\max_{0 \le \alpha < \beta \le \pi}
\left|S_p(\lambda;\alpha,\beta)\right| \ll  p^{3/4} (\log p)^{1/2} .
$$
\end{lemma}

\begin{proof} 
We  recall that for any real smooth function $F(\gamma)$ defined on 
the interval $[0,1]$ 
and any  sequence of
$N$ real numbers $\gamma_1,\ldots,\gamma_N\in[0,1]$ of discrepancy $D$, 
we have
$$
 \frac{1}{N} \sum_{n =1}^N F(\gamma_n) = \int_0^1 F(\gamma) d\, \gamma + O(D\max_{0\le \gamma  \le
1}\left|F'(\gamma) \right|),
$$
see~\cite[Chapter~2, Theorem~5.4]{KuNi}. 
Writing 
$$
S_p(\lambda;\alpha,\beta) =  \sum_{a \in \cA_p(\alpha,\beta)} \cos\(2\pi \frac{\lambda a}{p}\)
+ i \sum_{a \in \cA_p(\alpha,\beta)} \sin\(2\pi \frac{\lambda a}{p}\)
$$
and applying Lemma~\ref{lem:Discr},  we obtain the desired bound. 
\end{proof}

\section{Sato--Tate Conjecture for Sum Sets}
\label{sec:main}

\begin{theorem}
\label{thm:Set W} For any two sets $\cU, \cV \subseteq \F_p^*$, we have
$$
\max_{0 \le \alpha < \beta \le \pi} 
\left|\# \cW_p(\cU, \cV;\alpha,\beta)   - \mu_{ST}(\alpha,\beta) \#\cU\# \cV  \right| \le 
 \sqrt{\# \cU \# \cV} p^{3//4}  (\log p)^{1/2}. 
$$
\end{theorem}

\begin{proof}
Using the identity
\begin{equation}
\label{eq:Ident}
\frac{1}{p}\sum_{\lambda \in \F_p } \ep(\lambda c)=
\left\{\begin{array}{ll}
1&\quad\text{if $c=0$,}\\
0&\quad\text{if $c \in \F_p^*$,}
\end{array}
\right.
\end{equation}
we write 
\begin{eqnarray*}
 \# \cW_p(\cU, \cV;\alpha,\beta) & = &\sum_{u \in \cU} \sum_{v \in \cV}
 \sum_{a \in \cA_p(\alpha,\beta)} \frac{1}{p}\sum_{\lambda
\in
\F_p }
\ep(\lambda(u + v -a) 
\\ & = &  \frac{1}{p}\sum_{\lambda
\in
\F_p } S_p(-\lambda;\alpha,\beta)  \sum_{u \in \cU}\ep(\lambda u) \sum_{v \in \cV}
\ep(\lambda v). 
\end{eqnarray*} 
Separating the term $\# \cA_p(\alpha,\beta)\# \cU \# \cV/p$ corresponding 
to $\lambda = 0$,  we derive 
\begin{equation}
\label{eq:W and R}
 \# \cW_p(\cU, \cV;\alpha,\beta) = \frac{\# \cA_p(\alpha,\beta)\# \cU \# \cV}{p} +  O\( R\),
\end{equation}
where
$$
R = \frac{1}{p} \sum_{\lambda \in \F_p^* }  
\left|S_p(-\lambda;\alpha,\beta)\right|
 \left| \sum_{u \in \cU}\ep(\lambda u) \right| 
\left| \sum_{v \in \cV} \ep(\lambda v)\right|.
$$
Using Lemma~\ref{lem:Exp Sum} and the the Cauchy inequality, 
we obtain  
\begin{eqnarray*}
R  & \le  & p^{-1/4} (\log p)^{1/2} \sum_{\lambda \in \F_p^* }   
 \left| \sum_{u \in \cU}\ep(\lambda u) \right| 
\left| \sum_{v \in \cV} \ep(\lambda v)\right|\\
& \le  & p^{-1/4} (\log p)^{1/2} \(\sum_{\lambda \in \F_p^* }   
 \left| \sum_{u \in \cU}\ep(\lambda u) \right|^2\)^{1/2}
 \(\sum_{\lambda \in \F_p^* }    
\left| \sum_{v \in \cV} \ep(\lambda v)\right|^2\)^{1/2}.
\end{eqnarray*} 
Furthermore, by~\eqref{eq:Ident} we see that
$$
\sum_{\lambda \in \F_p^* } \left| \sum_{u \in \cU}\ep(\lambda u) \right|^2
\le  \sum_{\lambda \in \F_p}   
 \left| \sum_{u \in \cU}\ep(\lambda u) \right|^2 =
\sum_{\lambda \in \F_p }   
  \sum_{u_1,u_2 \in \cU}\ep(\lambda (u_1 - u_2)   = p \# \cU.
$$ 
Similarly, 
$$
\sum_{\lambda \in \F_p^* } \left| \sum_{v \in \cV}\ep(\lambda v) \right|^2\le p \# \cV.
$$
Collecting the above estimates together, 
we obtain 
$$
R \le \sqrt{\# \cU \# \cV}  p^{3/4} (\log p)^{1/2}
$$
which after substitution in~\eqref{eq:W and R} and using~\eqref{eq:Nied Bound}
leads us to the bound 
\begin{eqnarray*}
\lefteqn{
\left|\# \cW_p(\cU, \cV;\alpha,\beta)   - \mu_{ST}(\alpha,\beta) \#\cU\# \cV  \right|}\\ 
& & \qquad \qquad \qquad \ll   \#\cU\# \cV  p^{-1/4}+  \sqrt{\# \cU \# \cV}  p^{3/4}
(\log p)^{1/2}.
\end{eqnarray*}
It remains to note that 
$$
 \#\cU\# \cV  p^{-1/4} \le   \sqrt{\# \cU \# \cV}  p^{3/4}
$$
thus the first term never dominates. 
\end{proof}

Clearly the asymptotic formula of Theorem~\ref{thm:Set W} is nontrivial 
under the condition~\eqref{eq:Threshold}.

\section{Nonlinearity}

 \begin{theorem}
\label{thm:Nonlin}
For the  nonlinearity  $N(f)$ of the Boolean function  $f(h) $, 
given by~\eqref{eq:Bool},  we have
$$
N(f) = 2^{n-1} \(1  + O\( 2^{-n/16} n^{1/8} \)\).
$$
\end{theorem}

\begin{proof} We estimate the Fourier coefficients ${\widehat f}(k)$ of $f$
by using the result  that for any integers $M$, $h_1$, $h_2$ with 
$0 \le M \le M+c_1 < M+ c_2 < 2^n$ we have 
$$
\sum_{b=0}^{M-1} (-1)^{f(b+c_1)+ f(b+c_2)}  \ll 
M^{2/3} p^{1/6} (\log p)^{1/3} + p^{1/2} \log p , 
$$
which is a combination of~\cite[Lemma~2.3]{FMRS} 
with a special case $r = 2$  of~\cite[Lemma~4.4]{FMRS}. 
In fact, the above bound can be simplified as 
\begin{equation}
\label{eq:Corr}
\sum_{b=0}^{M-1} (-1)^{f(b+c_1)+ f(b+c_2)}  \ll M^{2/3} p^{1/6} (\log p)^{1/3}  
\end{equation}
(since for $M \le  p^{1/2} \log p $ the bound~\eqref{eq:Corr} 
is trivial and for $M >  p^{1/2} \log p $ we also have 
$ M^{2/3} p^{1/6} (\log p)^{1/3} >  p^{1/2} \log p$). 
  
 We now fix some $m \le n$  and write  $a,r \in \Bn$ as 
$$
a = b + 2^m  c  \qquad \text{and} \qquad r = s + 2^m t, 
$$
with $0 \le b,s< 2^m$ and  $0 \le c,t < 2^{n-m}$. In particular
$$
\langle a,r  \rangle = \langle b,s\rangle + \langle c,t \rangle.
$$
Therefore, 
\begin{eqnarray*}
|{\widehat f}(r)|  & =  &|{\widehat f}(s + 2^m t )| = \left|2^{-n} \sum_{b = 0}^{2^{m} -1} \sum_{c =
0}^{2^{n - m} -1}  (-1)^{ f( b + 2^m  c ) + \langle b,s \rangle + \langle c,t  \rangle} 
\right|\\
 & \le  &2^{-n}  \sum_{b = 0}^{2^{m} -1} 
\left| \sum_{c= 0}^{2^{n - m} -1}  (-1)^{  f( b + 2^m  c ) +\langle c,t \rangle }  \right|. 
\end{eqnarray*}
By the Cauchy inequality, we obtain 
\begin{eqnarray*}
|{\widehat f}(r)|^2 
& \le  &2^{m-2n }  \sum_{b = 0}^{2^{m} -1} 
\left| \sum_{j = 0}^{2^{n - m} -1}   (-1)^{  f( b + 2^m  c ) +\langle c,t \rangle }  \right|^2\\
& =  &2^{m-2n }  \sum_{b = 0}^{2^{m} -1} 
\sum_{c_1, c_2 = 0}^{2^{n - m} -1}   (-1)^{ f( b + 2^m  c_1 )+    f(b + 2^m  c_2) 
+ \langle c_1, t \rangle + \langle c_2,t \rangle }  \\
& \le &   
2^{m-2n }\sum_{c_1, c_2 = 0}^{2^{n - m} -1}    \left| \sum_{b = 0}^{2^{m} -1}  
 (-1)^{ f( b + 2^m  c _1)+    f(b + 2^m  c_2) }\right|.
\end{eqnarray*}
For $2^{n-m}$ choices of $c_1= c_2$ the sums over $b$ is equal to $2^m$.
For the other choices  of $c_1$ and $c_2$ we can use the bound~\eqref{eq:Corr}, 
getting 
\begin{eqnarray*}
|{\widehat f}(r)| ^2& =& O\( 2^{m-2n }\(2^{n-m} 2^m + 2^{2(n-m)} 2^{2m/3} 2^{n/6} n^{1/3}\)\)\\
& =& O\( 2^{m-n}+  2^{n/6 - m/3} n^{1/3}\).
\end{eqnarray*} 
We now define $m$ by the inequalities 
$2^{m}  \le  2^{7n/8} n^{1/4}  < 2^{m+ 1} $, 
and after simple calculations  conclude the proof. 
  \end{proof}

\section{Comments}

It seems very plausible that~\cite[Corollary~2.4]{Mich2} can be used to derive 
a nontrivial estimate for sums 
$$
T_p(\chi;\alpha,\beta) = \sum_{a \in \cA_p(\alpha,\beta)} \chi(a),
$$
with a nonprincipal multiplicative character $\chi$ of $\F_p^*$. 
In this case one can obtain a multiplicative analogue
of our results and study the  set
$$
\cZ_p(\cU, \cV;\alpha,\beta) = \{(u,v) \in \cU\times \cV\ : \ uv
\in \cA_p(\alpha,\beta)\}. 
$$
Multidimensional analogues of our results which involve joint distributions
of Kloosterman sums can be obtained as well. 

Also, as curiosity, we mention that Theorem~\ref{thm:Set W}
can be combined with the technique of~\cite{BaSh1,BaSh2,BaSh3} to study
sets of  elements of \emph{Beatty} sequence 
$\fl{\vartheta m + \rho}$ (where  $\vartheta > 0$ and $\rho$ are real) which belong to
$\cA_p(\alpha,\beta)$, that is, sets  of the form 
$$
\cB_p(\vartheta, \rho, M;\alpha,\beta) = \{m \in [1, M] \ : \ 
\fl{\vartheta m + \rho} \in
\cA_p(\alpha,\beta)\} .
$$

\end{document}